\documentclass[11pt,leqno]{article}
\usepackage{mathrsfs}
\usepackage{amsfonts}
\usepackage{amsmath}
\usepackage{amsthm}
\usepackage{amstext}
\usepackage{amsopn}
\usepackage[pdftex]{graphicx}
\usepackage{setspace}
\usepackage{algpseudocode}
\DeclareGraphicsExtensions{.eps, pdf, .jpg, .tif}

\newtheorem{theorem}{Theorem}[section]

\newtheorem{proposition}[theorem]{Proposition}

\begin{document}
\openup .32\jot
\title{Mathematical aspects of the combinatorial game ``Mahjong''}

\author{Yuan Cheng, Chi-Kwong Li\\
Department of Mathematics, College of William and Mary, \\ 
Williamsburg, VA 23187, USA. \\
ycheng04@email.wm.edu,  ckli@math.wm.edu, \\
Sharon H. Li \\
Microsoft Corporation, Redmond, WA 98052, USA.\\ 
lisharon@outlook.com
}

\date{}
\maketitle

\begin{abstract}
We illustrate how one can use basic combinatorial 
theory and computer programming technique (Python) to 
analyze the combinatorial game: Mahjong.
The results confirm some folklore concerning the game, 
and  expose some unexpected results. Related results
and possible future research in connection to 
artificial intelligence are mentioned. Readers interested in the subject
may  further develop the techniques to deepen the study
of the game, or study other combinatorial games.
\end{abstract}



\section{Introduction}

Mahjong is a popular recreational game which originated in China a long time ago 
\footnote{There are other theories of the origin of Mahjong. Many believed it was 
introduced 150 years ago but some say it was invented by Confucius 2500 years ago. 
\cite{re2,mahjong}}. Nowadays, it is widely played in different countries, including 
the United States.
It is a game of skill, strategy, calculation, and some luck. Some researchers have 
even
suggested that Mahjong is a good cognitive game with positive impact for patients 
with Alzhemier's disease; see \cite{Chengetal}.

The purpose of this note is to explore some mathematical aspects of the 
Mahjong game  using  elementary combinatorial theory and some 
basic programming technique  (Python).
Our study leads to  
affirmative answers on some folklore
concerning the game, and some unexpected 
results. This study is our first attempt to study the Mahjong game using 
mathematical and computational techniques. We will conclude the paper with some
possible future research in connection to artificial intelligence, and indicate the
difference between the Mahjong game, and other recreational games such as chess and the game of Go. Readers interested in the subject
can further develop the techniques to deepen the study
of the game, or study other combinatorial games.

\section{Basic rules and mathematics questions of Mahjong game}
\label{sec:1}
The rules of Mahjong are simple.\footnote{There are 
different variations of the game \cite{re3}. 
Here, we describe the basic version.}  There are 144 tiles in total, 
consisting 
of 36 tiles of bamboo type, 36 tiles of dot type, 36 tiles of character type, and some 
special tiles. There are 4 copies of each of the distinct bamboo, dot, character, and special tiles. In addition, there is 1 copy each of 4 flower and 4 season tiles. These tiles are shown in the image below:\footnote{Image from \cite{mahjong}}
\\ \includegraphics[width=\textwidth]{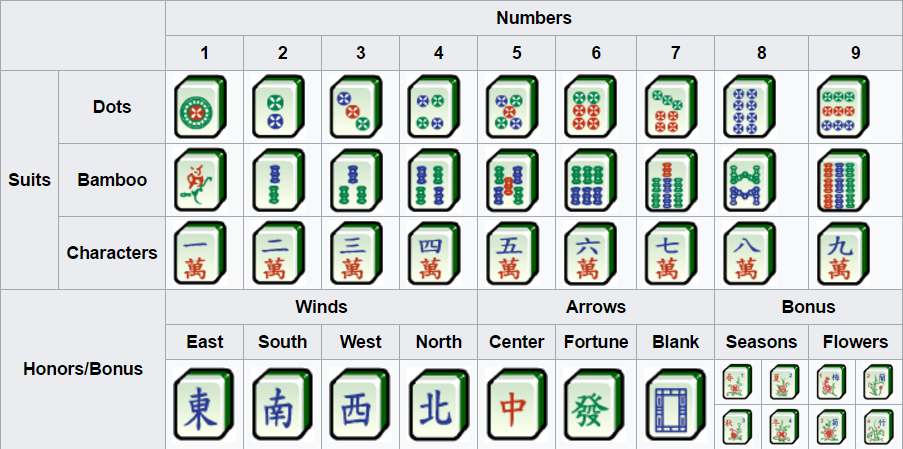}

When the game starts, each of the four players draws 13 tiles as a starting hand. Then, 
each player draws and then discards one tile in turns until one player forms a winning hand 
by using 13 tiles on hand and a newly drawn tile or a newly discarded tile of another player. 
A standard winning hand consists of an identical pair, and four sets of pungs or chows, where a pung is three identical tiles and a chow is three 
consecutive tiles from the same suit of dots, bamboo, or
characters.  The following picture 
shows two examples of a winning hand.\footnote{Images from \cite{re3}}
\\ \includegraphics[width=\textwidth]{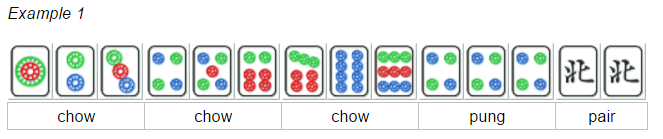} 
\\ \includegraphics[width=\textwidth]{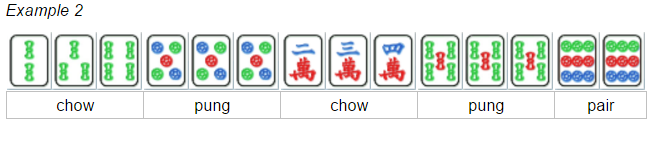} 
\\ The flowers and seasons tiles can add additional points to the winning score. Whenever one
 draws a season or flower tile, one puts it face up and draw another tile. 
Different winning hands will determine different winning scores. 
The score of the winner depends on how many seasons and 
flowers the player has and the rareness of the winning hand. 
For more details, one may see \cite{re3}.

A motivation of our study is the following hand:
\\ \includegraphics[width=\textwidth]{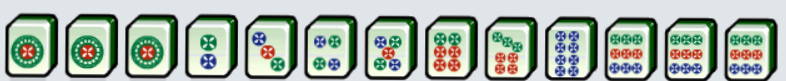}
\\ We will express this special hand as 
$X_1X_1X_1X_2X_3X_4X_5X_6X_7X_8X_9X_9X_9$ 
for  notational convenience. This hand is special because any 
additional dot tile would lead to a winning hand. For example, if we draw a one dot tile,
then we get a winning hand:
$$X_1X_1X_1 \quad X_1X_2X_3 \quad X_4X_5X_6 \quad X_7X_8X_9 \quad X_9X_9;$$
if we draw a dot 2, then we get a wining hand: 
$$X_1X_1X_1 \quad X_2X_2 \quad X_3X_4X_5 \quad X_6X_7X_8 \quad X_9X_9X_9.$$
One can check that each of the nine dot tiles can make this hand a winning hand.
This hand is called the ``Nine Gates".

It is believed that the ``Nine Gates" is unique. In other words,
if one has a hand of 13 tiles of dots, this is the only hand of dot tiles that 
yields a winning hand for any addition of a dot tile. However, there is no 
known mathematical proof of this folklore.

In connection to this, one can also get an ``Eight Gates" hand of dot tiles
a winning hand can be formed using eight of the nine dot tiles.
It is unclear how many such hands are possible, and which is the exceptional 
dot tile which cannot form a winning hand when added to the ``Eight Gates" hand. 
For instance, it is believed that there is no  
``Eight Gates'' hand so that one can win with any dot tile but the 5 dot tile. 

Of course, one can ask similar questions for ``Seven Gates", ``Six Gates'', etc.

In the next two sections, we will use combinatorial theory and 
develop a Python program to answer these questions and explore related results.
In Section 5, we describe the computational results and their implications;
related problems and further research will be discussed in Section 6.

One can download the Python program from 
http://cklixx.people.wm.edu/mathlib/Mahjong.py.
The computational results is contained in 
http://cklixx.people.wm.edu/mathlib/Mahjong-results.txt. 

\section{Mathematical Analysis}

In this section, we focus on Mahjong hands of 13 tiles chosen from the 
36 dot tiles to study the questions of ``Nine Gates'', ``Eight Gates'', etc.
We will continue to use the notation
$$X_1, \dots, X_9$$ 
to represent the 1-dot, \dots, 9-dot tiles each with 4 copies, 
and denote a hand by a  ``product'' of 
13 terms such as

\medskip\centerline
{$X_1X_1X_1X_2X_3X_4X_5X_6X_7X_8X_9X_9X_9$,  which may further simplify to
 $X_1^3X_2X_3X_4X_5X_7X_8X_9^3.$}

\medskip
We have the following facts about these 36 tiles based on 
basic combinatorial theory; for example, see \cite{re4}.

\begin{proposition} Consider 36 dot tiles with 4 copies of $X_1, \dots, X_9$.
Suppose a 13-dot hand is represented as $X_1^{n_1} \cdots X_9^{n_9}$
and associated with a sequence $(n_1, \dots, n_9)$ 
with $0 \le n_j \le 4$ for all $j$ such that $n_1 + \cdots + n_9 = 13$.

\begin{itemize}
\item[{\rm (a)}]  The number of ways to choose 13 random tiles from 36 tiles 
(allowing repeated patterns): $${36 \choose 13} = 2310789600.$$

\item[{\rm (b)}]  The number of ways of getting a certain 13 dot tiles hand 
with $n_j$ copies of $j$-dot tiles for $j = 1, \dots, 9$, so that $n_1, \dots, n_9$
are integers in $\{0,1,2,3,4\}$
adding up to 13: $${4 \choose n_1} \cdots {4 \choose n_9}.$$

\item[{\rm (c)}] The probability of getting a 13 dot tiles hand with $n_j$ copies of $j$-dot tiles:
$$\frac{{4 \choose n_1} \cdots {4 \choose n_9}}{{36 \choose 13}}.$$
So, the probability of getting the hand of nine gates
$X_1X_1X_1X_2 X_3X_4X_5X_6X_7X_8X_9X_9X_9$ out of the 36 dot tiles is 
$$\frac{262144}{2310789600}=0.00011344347.$$

\item[{\rm (d)}] All possible selections of $m$ tiles out of the 36 tiles for 
$m = 0, \dots, 36$, correspond to the degree $m$ terms in 
the expansion:
$$
(1+X_1+X_1^2+X_1^3+X_1^4)(1+X_2+X_2^2+X_2^3+X_2^4)\cdots(1+X_9+X_9^2+X_9^3+X_9^4)$$
$$ = \sum_{0 \le n_1 + \cdots + n_9 \le 36}  X_1^{n_1} X_2^{n_2} \dots X_9^{n_9}
\hskip 3in \ ~$$
In particular, the term with lowest degree
$1 = X_1^0 \cdots X_9^0$ corresponds to the selection of none of the
tiles $X_1 \dots X_9$, and the term with highest degree $X_1^4 \cdots X_9^4$
corresponds  to the selection of all the 36 tiles.

\item[{\rm (e)}] The number of different 13-dot hands 
equals to the total number of summands $X_1^{n_1}X_2^{n_2}\cdots X_9^{n_9}$ 
with $n_1+n_2+...+n_9=13$ in the expansion in {\rm (d)}, and equal to 
coefficient of $X^{13}$ in the expansion: 
$$(1+X+X^2+X^3+X^4)^9 = \sum_{i=0}^{36} \alpha_i X^i.$$
We have $\alpha_{13} = 93600.$

\item[{\rm (f)}] Define the dual of $X_{j_1}X_{j_2} \cdots X_{j_r}$ as
$X_{10-j_1} X_{10-j_2} \cdots X_{10-j_r}$.
It is easy to see that the dual of a chow, a pung or a pair is also a chow,
a pung or a pair, respectively. Consequently, the adding of 
$X_j$ to the 13-dot hand 
$X_1^{n_1} X_2^{n_2}\cdots X_9^{n_9}$ form a winning hand 
if and only if adding of $X_{10-j}$ form a winning hand of its dual.
\end{itemize}
\end{proposition}

The properties (a) -- (e) allow us to set up the Python program to determine the 
hands of nine gates, eight gates, etc. and compute their probability.
Property (f) is of theoretical interest that
every single hand of 13 tiles has a dual hand if we replace 
$X_k$ by $X_{10-k}$.
If the original hand can win with an $\ell$-dot tile, 
then its dual hand can win with $(10-\ell)$-dot tile. 
So two hands that are dual to each other can win by same number of tiles. 
For example, $X_1X_1X_1X_1 X_2X_2X_2X_2X_3X_3X_3X_4X_5$ and 
$X_5X_6X_7X_7X_7X_8X_8X_8X_8X_9X_9X_9X_9$ are dual to each other.
Adding $X_3, X_4$ or $X_6$ to the first hand will yield a winning hand.
Accordingly, adding $X_7$, $X_6$ or $X_4$ to the second hand will yield a winning hand.
Evidently, the dual hand of the ``Nine Gates" is itself.

\section{Programming}

Based on the mathematical results in the previous section, we write a Python program to 
study different hands of 13 dot tiles. The program is available at

http://cklixx.people.wm.edu/mathlib/Mahjong.py.

\medskip\noindent
One can also see the listing of the results at

http://cklixx.people.wm.edu/mathlib/Mahjong-results.txt.

The basic idea of this program is to generate all $93600$ of such hands
calculated in Proposition 3.1 (e), 
and test each of them to see how many different tiles would complete a winning hand.
It is worth pointing that in our Algorithm 1, we associate each hand of 13 
dot tiles $X_1^{n_1} \cdots X_9^{n_9}$
with the sequence $(n_1, \dots, n_9)$ such that $0 \le n_j \le 4$ for every
$j$ and  $n_1 + \cdots + n_{9} = 13$. This allows us to modify the program 
easily to check what are needed to form a winning pattern for a reduced hand
after some ``pungs'' or ``chows'' were performed in a game. 

For each 13-tile hand, we add a new tile from 1-9 to it, to create a 14-tile hand.  
To determine if this is a winning hand, we must identify a pair and four sets of pungs and 
chows. For each of the tiles that appeared at least twice in the hand, we take two of them out 
as the pair. If the remaining 12-tile hand $\{j_1, \dots, j_{12}\}$
can be divided into four sets of pungs and chows, then this is a winning hand.
The following proposition is useful for the test.

\begin{proposition}
Use the notation in the previous discussion.
If $j_1 = j_2 = j_3$, we may always assume that they form a pung and check whether 
the remaining pieces $\{j_4, \dots, j_{12}\}$ form three sets of pungs and chows.
\end{proposition}

\it Proof. \rm
Assume that $j_1 = j_2 = j_3$, and we can divide $\{j_1, \dots, j_{12}\}$
into three sets of pungs and chows without using $\{j_1,j_2,j_3\}$ as a pung.
Then $j_1, j_2, j_3$ will be in three sets of chows
of the form $\{j_1, j_1+1, j_2 + 2\}$ plus one other set of pung or chow.
But then then we can rearrange the twelve tiles as three sets of chows 
into three sets of pungs $\{j_1,j_1, j_1\}, \{j_1+1, j_1+1,j_1+1\}, \{j_1+2, j_1+2, j_1+2\}$
 together with the remaining set of pung or chow.  Thus, our claim is proved. \qed
 
\medskip
Clearly, after removing $\{j_1, j_2, j_3\}$ from $\{j_1, \dots, j_{12}\}$,
if the three smallest number in the remaining set are the same, we can again
assume that they form a pung. Else, we will extract a set of chow and proceed
in a similar manner. We will use this idea in Algorithm 3 below.

\medskip
 Below is the pseudocode of our algorithms:

\noindent{\bf Algorithm 1: Pseudocode for generating every possible hand.}
\\ Goal: Form every possible hand of 13 tiles.
\\ Output: All possible 13-tile hands
\begin{algorithmic}[1]
\State $hand \gets$ array of size 9 
\Comment{$hand[i]$ represents number of i-dot tiles, 0-4}
\For {each hand where hand[i] from 0-4}
    \If {$sum(hand) = 13$}
    \State add $hand$ to $allPossibleHands$
    \EndIf
\EndFor
\State\Return $allPossibleHands$
\end{algorithmic}
\medskip

\noindent{\bf Algorithm 2: Pseudocode for checking for number of gates.}
\\ Goal: Check whether a 13-tile hand can form a winning pattern with 1 more tile.
\\ Input: A 9-element list called $hand$, representing 13 tiles.
\\ Output: Which tiles does $hand$ need to win.
\begin{algorithmic}[1]
\State $isWinningTile \gets$ array of size 9 initialized to false
\For {$i = 1,\dots,9$}
    \If { $hand[i] \neq 4$}
    	\State $hand[i] \gets hand[i] + 1$ \Comment{Add an i-tile to get a 14-tile hand}
    	\For {$j = 1,\dots,9$, if $hand[j] \geq 2,$}
    		\State $hand[j] \gets hand[j] -2 $ \Comment{first select $j$ as a pair}
    		\If { $checkFourSets(hand)$} \Comment{check if other 12 tiles has 4 sets}
    			\State $isWinningTile[i] \gets$ True
    		\EndIf
    	\EndFor
    \EndIf
\EndFor
\State\Return $isWinningTile$
\end{algorithmic}
\medskip

\noindent{\bf Algorithm 3: Pseudocode for checking four sets of pungs or chows}
\\ Goal: Check whether a 12-tile hand consists of four sets of 3-tile.
\\ Input: A 9-element list called $hand$, representing 12 tiles.
\\ Output: Whether the hand forms four sets of pungs or chows
\begin{algorithmic}[1]
\State $setsFound = 0$
\For {$i = 1,\dots,9$}
    \If { $hand[i] \geq 3$} \Comment{Check for pung}
    	\State $hand[i] \gets hand[i] - 3$
    	\State $setsFound = setsFound + 1$
    \EndIf
    \If { $i+2 < len(hand)$}  \Comment{Check for chow}
    	\State $minThree = min(hand[i], hand[i+1], hand[i+2])$
    	\State $hand[i] = hand[i] - minThree$
    	\State $hand[i+1] = hand[i+1] - minThree$
    	\State $hand[i+2] = hand[i+2] - minThree$
    	\State $setsFound = setsFound + minThree$
    \EndIf
\EndFor
\If { $setsFound =4$}  
\State\Return $True$
\Else
\State\Return $False$
\EndIf
\end{algorithmic}

\section{Computational Results and Implications}
The program gives us the following results.

\begin{proposition}
Consider the 93600 different 13 dot hands.
\begin{itemize}
\item[{\rm (a)}] 
The ``Nine Gates" $X_1X_1X_1X_2X_3X_4X_5X_6X_7X_8X_9X_9X_9$
 is the unique hand winning all 9 pieces 
 with the probability of 0.000113  for a 13 dot hands as 
 shown in Proposition 3.1 (c). 

\item[{\rm (b)}] 
There are 16 hands winning 8 pieces with a combined probability 
0.0001 
of drawing.{\rm
\\ $X_3X_3X_3X_4X_5X_6X_7X_8X_8X_8X_9X_9X_9$ [winning except for the 1 dot tile]
\\ $X_3X_3X_3X_4X_5X_5X_6X_6X_7X_7X_8X_8X_8$ [winning except for the 1 dot tile]
\\ $X_3X_3X_3X_4X_4X_5X_5X_6X_6X_7X_8X_8X_8$ [winning except for the 1 dot tile]
\\ $X_2X_3X_4X_4X_4X_4X_5X_6X_7X_8X_9X_9X_9$ [winning except for the 4 dot tile]
\\ $X_2X_3X_3X_3X_3X_4X_4X_5X_6X_7X_8X_8X_8$ [winning except for the 3 dot tile]
\\ $X_2X_2X_2X_3X_4X_5X_6X_7X_7X_7X_9X_9X_9$ [winning except for the 9 dot tile]
\\ $X_2X_2X_2X_3X_4X_5X_6X_7X_7X_7X_8X_8X_8$ [winning except for the 9 dot tile]
\\ $X_2X_2X_2X_3X_4X_5X_6X_7X_7X_7X_7X_8X_9$ [winning except for the 7 dot tile]
\\ $X_2X_2X_2X_3X_4X_5X_6X_6X_7X_7X_7X_7X_8$ [winning except for the 7 dot tile]
\\ $X_2X_2X_2X_3X_4X_4X_5X_5X_6X_6X_7X_7X_7$ [winning except for the 9 dot tile]
\\ $X_2X_2X_2X_3X_3X_4X_4X_5X_5X_6X_7X_7X_7$ [winning except for the 9 dot tile]
\\ $X_2X_2X_2X_3X_3X_3X_4X_5X_6X_7X_8X_8X_8$ [winning except for the 1 dot tile]
\\ $X_1X_2X_3X_3X_3X_3X_4X_5X_6X_7X_8X_8X_8$ [winning except for the 3 dot tile]
\\ $X_1X_1X_1X_3X_3X_3X_4X_5X_6X_7X_8X_8X_8$ [winning except for the 1 dot tile]
\\ $X_1X_1X_1X_2X_3X_4X_5X_6X_6X_6X_6X_7X_8$ [winning except for the 6 dot tile]
\\ $X_1X_1X_1X_2X_2X_2X_3X_4X_5X_6X_7X_7X_7$ [winning except for the 9 dot tile]}

\item[{\rm (c)}] There are 79 hands of ``Seven Gates'' with a combined probability
0.000942  of drawing. 

\item[{\rm (d)}] There are 392 hands of ``Six Gates''
with a combined probability 0.005408 of drawing.

\item There are 1335 hands of  ``Five Gates'' with a combined probability
0.014215 of drawing.

\item There are 2948 hands of ``Four Gates''  
with a combined probability 0.029812 of drawing. 

\item There are 6739 hands of ``Three Gates'' with a combined probability 
0.097559  of drawing.

\item There are 14493 hands of ``Two Gates''
with a combined probability 0.178968  of drawing.

\item There are 14067 hands of  ``One Gate'' 
with a combined probability 0.148473  of drawing.

\item There are 53530 hands which cannot win
with any additional piece, with a combined 
probability 0.524409 of drawing.
\end{itemize}
\end{proposition}

\medskip
Several remarks are in order in connection to Proposition 4.1.
\begin{enumerate}
\item It is confirmed in (a) that the  ``Nine Gates'' is the unique 
hand of 13 dots that can form a winning hand with the addition of any dot tile,
 and the probability
of getting such a hand is  
$$4^9/{36 \choose 13}=0.000113.$$

\item In (b), one can compute the numbers of ways to get each of the 
hand of 13 dot tiles using Proposition 3.1 (b).

For example, 
the number of ways to get $X_3X_3X_3X_4X_5X_6X_7X_8X_8X_8X_9X_9X_9$ equals
$${4\choose 0}^2{4\choose 3}{4\choose 1}^4{4\choose 3}^2 = 4^7,$$
and the number of ways to get $X_3X_3X_3X_4X_5X_5X_6X_6X_7X_7X_8X_8X_8$ equals
$${4\choose 0}^2{4\choose 3}{4\choose 1}{4\choose 2}^3{4\choose 1}
{4\choose 0} = 4^36^3.$$
Adding the number of ways to get the 16  ``Eight Gates'' hands, we have
$$4^7+4^36^3 + 4^36^3 + 4^7 + 4^5 6 + 4^7 + 4^7 + 4^7 + 4^56 + 4^36^3
+ 4^36^3 + 4^7 + 4^7 + 4^7 + 4^7 +4^7 = 231424.
$$
Dividing the sum by ${36 \choose 13}$,
we see that the probability of $0.0001$ of drawing these hands.

\item
It is somewhat interesting (and counter intuitive) 
that if one draws 13 tiles out of the 36 dot tiles,
there is a higher probability of getting the  ``Nine Gates'' 
(0.000113) 
than that of getting one of the ``Eight Gates'' (0.0001).

\item
Note that the 8 of the 16 hands in (b) are dual to
the other 8 hands as described in Proposition 3.1 (f).

\item 
It is also interesting to note that the only way to get an  ``Eight Gates'' hand
so that the 4-dot tile cannot be added to form a winning hand is:
$X_2X_3X_4X_4X_4X_4X_5X_6X_7X_8X_9X_9X_9$, where all the 4-dot tiles are used.

\item Similar comment applies to the hand where the $k$-dot tile 
cannot be added to form a winning hand for $k = 3,6,7$. In particular, 
when $k=3$ there are two such  ``Eight Gates''
hand. The same is true for $k = 7$.

\item There are three tiles that every ``Eight Gates" can win with. They are $2,5,8$.

\item By (b), we see that there is no  ``Eight Gates'' hand that win every
dot tiles except the $k$-dot for $k = 2, 5$ or $8$.

\item There are too many hands corresponding to ``Seven Gates'', ``Six Gates'', 
``Five Gates'', etc.
We do not list them in the proposition. Nevertheless, in the 
results output from our program in the Appendix, we put the statistics of
the number of hands corresponding to each of the n-Gates and indicate the  
tiles whose addition will form winning hands.
One can run the Python program available at 
http://cklixx.people.wm.edu/mathlib/Mahjong.py. 
to see all the possible outcomes as shown in 
http://cklixx.people.wm.edu/mathlib/Mahjong-results.txt.

\item For the five gates hand, there are 8 hands 
such that $1,5,9$ cannot be the winning 
pieces, namely,
$$
X_3X_4X_5X_5X_5X_5X_6X_7X_7X_7X_9X_9X_9, \quad
X_3X_4X_5X_5X_5X_5X_6X_7X_7X_7X_8X_8X_8, $$
$$
X_3X_3X_3X_4X_5X_5X_5X_5X_6X_7X_7X_8X_9, \quad
X_3X_3X_3X_4X_5X_5X_5X_5X_6X_7X_9X_9X_9, $$
$$
X_2X_2X_2X_3X_3X_3X_4X_5X_5X_5X_5X_6X_7, \quad
X_1X_2X_3X_3X_4X_5X_5X_5X_5X_6X_7X_7X_7, $$
$$
X_1X_1X_1X_3X_4X_5X_5X_5X_5X_6X_7X_7X_7, \quad
X_1X_1X_1X_3X_3X_3X_4X_5X_5X_5X_5X_6X_7.$$
In each case, $X_5$ appears four times so that $X_5$ cannot be the additional
tile to form a winning hand.

\item In the results output from our program in the Appendix, we also list
the statistics for ``Four Gates'' and ``Three Gates''.
Here we list the tiles whose addition will lead to a winning hand.

\item There are many known ``Three Gates'' hands. For every triple $(i,j,k)$ 
with $1 \le i < j < k \le 9$, one may ask whether there is a hand of 13 dot tiles
such that one can get a winning hand by adding a tile from $\{X_i, X_j, X_k\}$.
Our results show that out of the ${9 \choose 3} = 84$ possible choices of 
$\{X_i,X_j,X_k\}$
one can get  ``Three Gates'' hands with these sets of winning tiles 
with the following 11 exceptions: 
$$\{X_1,X_2,X_9\},\{X_1,X_3,X_8\},
\{X_1,X_5,X_7\},\{X_1,X_5,X_9\},
\{X_1,X_6,X_8\},\{X_1,X_8,X_9\},$$
$$
\{X_2,X_4,X_8\},\{X_2,X_4,X_9\},
\{X_2,X_6,X_8\},\{X_2,X_7,X_9\},
\{X_3,X_5,X_9\}.$$

\item It is easy to see that for any $1 \le i, j \le 9$, one can 
create a hand of 13 dot tiles so that a winning hand will  be formed by
adding the $i$-dot or the $j$-dot tile. 
For example, one may have a pair of $i$-dot tiles,
and a pair of $j$-dot tiles, and three sets of ``pung'' so that only the addition
of the $i$-dot or $j$-dot will lead to a winning hand.

\item Similarly, one can have a hand of 4 sets of  ``pung''  with a single $k$-dot tile
so that one can only use an additional $k$-dot to form a winning hand.
\end{enumerate}

\section{Related problems and Further Research}

We can also use our program to answer other problems. For example,
if one randomly draws 14 tiles from the 36 dot tiles, what is the probability
of getting a winning hand?

A calculation in our computer program shows that there are 118800 possible
14 tile dot hands, of which 13259 are winning.
As a result, the probability of getting 14 tiles that form a winning hand is: 
$0.11161$, which is larger than $\frac{1}{9}$. 
This result is higher than many Mahjong players would expect.

Of course, one can consider the full set of Mahjong with 144 tiles. 
Computing the probability of getting certain special hands will be
more complicated. 

In fact, a more challenging project is to develop a computing Mahjong-playing 
system. There has been great progress in research in  
artificial intelligence and machine learning.
Computer systems have been built that  beat the best chess player and Go player;
see \cite{Nielsen,Silver,deepblue,alphago}. It would be interesting to develop a 
Mahjong playing machine to beat the best Mahjong player in the world. 

Note that Mahjong is different from chess and Go because the players do not have
the full information of other players during the game.
One needs to anticipate what other players are hiding in their hands and 
create their own game plan. Also, skillful players would be able 
to anticipate other players' strategies by observing their discarded tiles. 
Also, it is possible for two or three players to form a coalition to play
against the other players. So, playing Mahjong well 
would require good use of combinatorial theory, probability theory,
game theory, psychology, etc.  To develop a good Mahjong 
playing machine will require another level of intelligence.

\medskip\noindent
{\bf Additional notes.}

\begin{itemize}
\item After this paper was submitted, the second author reported
the findings of the paper in Taiwan in a lecture tour
in 2018; he was alerted that
Taiwanese have different rules for Mahjong. For instance, 
a basic winning pattern requires 17 tiles consisting of 
one more set of chow or pung to the type described in this paper.
For example,
$$X_1X_1X_1X_2X_2X_2X_3X_3X_3X_4X_5X_6X_7X_8X_8X_8$$
together with any $X_j$ for $j = 1, \dots, 9$,
will form a winning hand.

\item One can modify our theory and Python program to find 
out all the 16-dot hands that will win with the addition of 
a tile from a certain set $\{X_{j_1}, \dots, X_{j_r}\}$.
The numerical results can be found at
http://cklixx.people.wm.edu/mathlib/Mahjong-results2.txt.

\item  It is interesting to point out that
there are 11 different hands of nine-gates in this case:

\qquad
[2, 2, 2, 3, 4, 5, 6, 7, 7, 7, 8, 8, 8, 9, 9, 9] 

\qquad 
[2, 2, 2, 3, 3, 3, 4, 4, 4, 5, 6, 7, 8, 9, 9, 9] 

\qquad
[1, 1, 2, 2, 2, 3, 3, 3, 4, 5, 6, 7, 8, 9, 9, 9] 

\qquad
[1, 1, 1, 2, 3, 4, 5, 6, 7, 7, 7, 8, 8, 8, 9, 9]

\qquad
[1, 1, 1, 2, 3, 4, 5, 6, 6, 7, 7, 8, 8, 9, 9, 9] 

\qquad
[1, 1, 1, 2, 3, 4, 5, 6, 6, 6, 7, 7, 7, 8, 8, 8]

\qquad
[1, 1, 1, 2, 3, 4, 5, 5, 6, 6, 7, 7, 8, 9, 9, 9] 

\qquad
[1, 1, 1, 2, 3, 4, 4, 5, 5, 6, 6, 7, 8, 9, 9, 9]

\qquad
[1, 1, 1, 2, 3, 3, 4, 4, 5, 5, 6, 7, 8, 9, 9, 9] 

\qquad
[1, 1, 1, 2, 2, 3, 3, 4, 4, 5, 6, 7, 8, 9, 9, 9]

\qquad
[1, 1, 1, 2, 2, 2, 3, 3, 3, 4, 5, 6, 7, 8, 8, 8]

\item
If one randomly picks 17 tiles from the 36 dot tiles,
there are 26414 out of the 175725 hands corresponding to 
winning hands. Thus, the probability is
0.15314. 

\item Professor Jun-Yi Guo of the National 
Taiwan Normal University used the idea of this paper together
with the theory of generating function to determine more 
general winning patterns of Mahjong.

\item The second author attended the SIAM  
(Society for Industrial and Applied Mathematics) Annual
Meeting at Portland, Oregon, July 9-13, 2018, 
and met two colleagues, Kon Aoki and Kochi Nakajima of Colorado
College. They gave a presentation on  
``Theoretical Mahjong:
Discard Piles Algorithm Simulator''.
The two colleagues pointed out that there are some Japanese
researchers doing mathematical
and A.I. study of Mahjong.

\end{itemize}

\medskip\noindent


%
%

\medskip\noindent
{\large\bf Acknowledgement}

The authors would like to thank the referee for the careful reading of the
manuscript and some helpful suggestions.
This research of the first two authors
was partially supported by the NSF EXTREEMS-QED grant 1331021,
and the Cissy Patterson Fund of the College of William and Mary.

%
%

\end{document}